\documentclass{amsart} 

\usepackage[english,french]{babel}
\usepackage[latin1]{inputenc}

\usepackage{amsmath}
\usepackage{amsthm}
\usepackage{amssymb}
\usepackage{tikz}

\usepackage{imakeidx}  
\usepackage{appendix}
\usepackage{tikz-cd}
\usepackage{verbatim}
\usepackage{enumitem}
\usepackage{multicol}

\usepackage[backref=page]{hyperref}
\usepackage{cleveref}

\hypersetup{colorlinks=true,linkcolor=blue,anchorcolor=blue,citecolor=blue}

\usepackage{adjustbox}

\usepackage{mathtools}

\newtheorem{thm}{Th\'eor\`eme}[section]

\newtheorem{conj}[thm]{Conjecture}

\newtheorem{prob}[thm]{Probl\`eme}

\theoremstyle{definition}

\theoremstyle{remark}

\numberwithin{equation}{section}

\newcommand{\R}{\mathbb R}
\newcommand{\Br}{{\rm Br}}
\newcommand{\Q}{\mathbb Q}
\newcommand{\F}{\mathbb F}
\newcommand{\C}{\mathbb C}
\newcommand{\Z}{\mathbb Z}

\renewcommand{\P}{\mathbb P}

\newcommand{\A}{\mathbb A}

\newcommand{\cl}{\overline}

\renewcommand{\phi}{\varphi}

\title[Une liste de probl\`emes]{Une liste de probl\`emes}

\author{Jean-Louis Colliot-Th\'el\`ene}

\address{Universit\'e Paris-Saclay, CNRS, Laboratoire de math\'ematiques d'Orsay, 91405, Orsay, France}
\email{jean-louis.colliot-thelene@universite-paris-saclay.fr}

 \date{10 d\'ecembre 2022}

\begin{document}
	\maketitle
	
	\section*{Introduction}
	
Dans cette note, je rassemble une liste de probl\`emes, la plupart bien connus, et rest\'es ouverts
 depuis de nombreuses ann\'ees. J'ai r\'efl\'echi \`a la plupart d'entre eux mais n'en revendique pas
 la propri\'et\'e.  Je mentionne certaines solutions partielles, sans faire un
 rapport syst\'ematique.  Je renvoie \`a \cite{CT87, CT98,  Sk01, CT03, CT11, SD11, CT19, W18}  pour cela.

 Les probl\`emes portent presque tous sur la g\'en\'eralisation en dimension plus grande que 1
 des deux \'enonc\'es suivants :
 
 Une conique lisse sur un corps $k$ qui poss\`ede un point rationnel sur $k$
 est isomorphe, sur $k$, \`a la droite projective $\P^1_{k}$. 
 Ceci donne une param\'etrisation biunivoque des points rationnels de la
 conique par les points rationnels de la droite projective.
 
 Sur un corps de nombres $k$, si une conique lisse a un point rationnel sur 
 tous les compl\'et\'es $k_{v}$ de $k$, alors elle a un point rationnel sur $k$.
Gr\^{a}ce \`a Hensel, ce crit\`ere est effectif. 
 
 Le premier \'enonc\'e remonte \`a l'Antiquit\'e, 
 le second fut \'etabli par Legendre
 sur les rationnels et par Hilbert sur les corps de nombres,
 et fut \'etendu par Minkowski et par Hasse aux quadriques
 de dimension quelconque.
 
 Beaucoup des probl\`emes mentionn\'es ici ont leur source dans mes travaux avec  Jean-Jacques Sansuc et avec 
Peter Swinnerton-Dyer dans les ann\'ees 1970 et 1980.
Certains des probl\`emes, en particulier ceux
sur les intersections de deux quadriques,  avaient fait l'objet
de rapports  non publi\'es en  1988 et  en 2005.

\section{Vari\'et\'es rationnelles et vari\'et\'es proches}

Soit $k$ un corps alg\'ebriquement clos.
On dit qu'une vari\'et\'e int\`egre $X$   sur $k$   est   {\it rationnelle} si elle est  birationnelle \`a un espace projectif $\P^d_{k}$, i.e. si son corps des fonctions $k(X)$ est transcendant pur sur $k$.

Parmi les exemples classiques de vari\'et\'es rationnelles, on trouve les vari\'et\'es sous-jacentes
\`a un groupe alg\'ebrique lin\'eaire connexe, et les vari\'et\'es projectives qui sont des espaces
homog\`enes de tels groupes.  Les quadriques lisses de dimension au moins~1 rentrent dans ce cadre.

 On dit qu'une vari\'et\'e   int\`egre $X$  sur $k$   est {\it unirationnelle} 
 s'il existe une application rationnelle dominante d'un espace projectif vers $X$.
 
En dimension 1 et en caract\'eristique z\'ero en dimension 2, unirationalit\'e
implique rationalit\'e. C'est faux d\`es la dimension 3 (Clemens--Griffiths,
 Iskovskikh--Manin, Artin--Mumford). 
 
 Parmi les exemples classiques de vari\'et\'es unirationnelles, on trouve 
  les quotients $G/H$ d'un groupe
lin\'eaire connexe $G$ par un sous-groupe ferm\'e  $H$ quelconque,
 non n\'ecessairement connexe. Pour $H$  fini, on conna\^{\i}t des exemples
 de tels quotients qui ne sont pas rationnels (Saltman, Bogomolov).
 
Dans la classification birationnelle des vari\'et\'es de dimension sup\'erieure d\'eve\-lop\-p\'ee vers 1990
(travaux de Koll\'ar, Miyaoka, Mori),  ce qui en dimension quelconque joue le r\^ole
des surfaces rationnelles dans la classification des surfaces, ce sont les
  vari\'et\'es rationnellement connexes \cite{K99, AK03}. 
  L'une des d\'efinitions, en caract\'eristique z\'ero, est que par deux points (ferm\'es) g\'en\'eraux d'une telle vari\'et\'e il passe
 une courbe de genre z\'ero, et ce sur tout corps alg\'ebriquement clos contenant 
 $k$.
  En caract\'eristique quelconque, la bonne d\'efinition est celle
 de vari\'et\'e s\'eparablement rationnellement connexe. Les deux notions co\"{\i}n\-cident en caract\'eristique nulle.
 Dans la suite de ce texte,  par  vari\'et\'e  rationnellement connexe on entendra vari\'et\'e 
 s\'eparablement rationnellement connexe. 
 
 Supposons $k$ de caract\'eristique z\'ero. Une vari\'et\'e unirationnelle est
 rationnellement connexe. La r\'eciproque est un grand probl\`eme ouvert.
Une vari\'et\'e  lisse, projective, lisse,  connexe, \`a fibr\'e anticanonique ample
est appel\'ee vari\'et\'e de Fano.
Un   th\'eor\`eme  important (Campana, Koll\'ar--Miyaoka--Mori) dit qu'une vari\'et\'e de Fano
  est rationnellement connexe.
  Ainsi toute hypersurface lisse $X \subset \P^n_{k}$, $n \geqslant 2$,  de degr\'e $d$ avec $d\leqslant n$ 
  est rationnellement connexe. Il en est donc ainsi des hypersurfaces cubiques lisses dans $\P^n_{k}, n \geqslant 3$.
  Il en est aussi ainsi des intersections compl\`etes lisses de deux quadriques dans $\P^n_{k}$, $n \geqslant 4$.
Ces derni\`eres sont des vari\'et\'es rationnelles sur le corps alg\'ebriquement clos $k$.
 
  Un autre th\'eor\`eme important   (Graber--Harris--Starr, 2003) dit que l'espace total
  d'une fibration de base rationnellement connexe et de fibres g\'en\'erales
   rationnellement connexes est  rationnellement connexe.
  
\medskip

  Soient maintenant $k$ un corps quelconque et $\cl{k}$ une cl\^{o}ture alg\'ebrique.
  Nous adoptons ici les conventions suivantes.
  
  On dit qu'une $k$-vari\'et\'e g\'eom\'etriquement int\`egre $X$
  est   rationnelle,  resp.   rationnellement connexe,
   si  la $\cl{k}$-vari\'et\'e $\cl{X}:=X \times_{k}\cl{k}$ est
  rationnelle,   resp.  rationnellement connexe.
  
  On dit qu'une $k$-vari\'et\'e g\'eom\'etriquement int\`egre $X$
  est $k$-rationnelle si $X$ est $k$-birationnelle \`a $\P^d_{k}$,
  i.e. le corps des fonctions $k(X)$  de $X$ est transcendant pur.
  
  On dit qu'une $k$-vari\'et\'e g\'eom\'etriquement int\`egre $X$
  est stablement $k$-rationnelle s'il existe des entiers $n \geqslant 0$ et $m \geqslant 0$
  tels que $X\times_{k}\P^n_{k}$ est $k$-birationnelle \`a $\P^m_{k}$.
    
  On dit qu'une $k$-vari\'et\'e g\'eom\'etriquement int\`egre $X$ de dimension $d$
  est    $k$-uni\-ration\-nelle s'il existe une application $k$-rationnelle
  dominante  d'un espace projectif $\P^d_{k}$ vers $X$.

\section{Principe de Hasse, approximation faible, obstruction de Brauer--Manin}

\'Etant donn\'ee une vari\'et\'e alg\'ebrique $X$ sur un corps $k$,
on note $X(k)$ l'ensemble de ses points rationnels. Pour $K/k$
une extension quelconque  de corps, on note $X(K)$ l'ensemble
des points rationnels sur $K$.

On veut donner des crit\`eres si possible effectifs permettant
de d\'ecider si une $k$-vari\'et\'e donn\'ee $X$ poss\`ede
un point rationnel.

Pour \'eviter des r\'ep\'etitions, on va d\'efinir un certain
nombre de propri\'et\'es.

\medskip

$({\bf PR}_{X})$  L'ensemble $X(k)$  est non vide, i.e. la $k$-vari\'et\'e poss\`ede un point 
 rationnel sur $k$.

\medskip

On s'int\'eresse particuli\`erement  au cas des corps finis, des
corps locaux (les corps $p$-adiques, les corps de s\'eries formelles
en une variable sur un corps fini, le corps $\R$ des r\'eels, le corps $\C$
des complexes), et des corps globaux (corps de nombres ou corps de fonctions
d'une variable sur un corps fini). Etant donn\'e un corps global $k$, et une place
$v$ de ce corps, on note $k_{v}$ le corps local compl\'et\'e par rapport \`a la place $v$.

Soit d\'esormais $k$ un corps global. 

  \`A toute  vari\'et\'e alg\'ebrique $X$ sur   $k$
on associe l'espace $X(\A_{k})$ de ses ad\`eles.
C'est un sous-ensemble du produit $\prod_{v}X(k_{v})$,
non vide si  ce produit est non vide.
 L'espace 
  $X(\A_{k})$ est muni d'une topologie naturelle. 
Si $X$ est projective, alors on a  $X(\A_{k})= \prod_{v}X(k_{v})$,
et la topologie de l'espace des ad\`eles co\"{\i}ncide avec la topologie
produit sur $\prod_{v} X(k_{v})$.

On introduit la propri\'et\'e :

\medskip

 $({\bf PH}_{X})$ Soit  on a $X(\A_{k} ) = \varnothing $, soit  on a $X(k) \neq \varnothing$.

  On dit que le principe de Hasse vaut pour une classe  
de vari\'et\'es alg\'ebriques d\'efinies sur  $k$  si, pour toute vari\'et\'e $X$ dans cette classe,
on a la propri\'et\'e  ${\bf PH}_{X}$.

\medskip

Pour $X/k$ lisse et g\'eom\'etriquement int\`egre, on introduit la propri\'et\'e d'approximation faible :

\medskip

$({\bf AF}_{X})$ L'image de l'application diagonale $X(k) \to  \prod_{v} X(k_{v})$ est dense.

\medskip

Cette propri\'et\'e implique ${\bf PH}_{X}$. Si $X(k)$ est non vide, elle \'equivaut au
fait que pour tout ensemble fini $S$ de places de $k$, l'ensemble $X(k)$
est dense dans le produit fini $\prod_{v\in S} X(k_{v})$. Il convient de noter
que pour certaines classes de vari\'et\'es,  on ne sait  pas  \'etablir  ${\bf PH}_{X}$
pour $X$ dans cette classe, mais que, sous l'hypoth\`ese $X(k)\neq \varnothing$,
la propri\'et\'e 
${\bf AF}_{X}$ est facile \`a \'etablir. 

\medskip

Pour $X/k$ non n\'ecessairement projective, on peut encore consid\'erer une variante de la propri\'et\'e
${\bf AF}_{X}$. Il s'agit du probl\`eme de l'approximation forte. Depuis 2008, il a \'et\'e 
\'etudi\'e  du point de vue de l'obstruction de Brauer--Manin,  dans plusieurs articles par Fei Xu et moi, Harari,
Borovoi, Demarche, Dasheng Wei, Yang Cao,  mais nous ne
le discuterons pas dans ce texte. Je renvoie \`a  \cite{BD13} et au rapport de Wittenberg \cite[\S 2.7, \S 3.2.4,  \S 3.3.4, \S 3.4.5]{W18} pour
des r\'esultats et r\'ef\'erences.
 
 \medskip
 
 Pour $X/k$ lisse et g\'eom\'etriquement int\`egre,  avec $X(k)\neq  \varnothing$,
 il y a lieu d'introduire la propri\'et\'e d'approximation ``faible faible''   \cite[Chap. 3]{Se92}):
 
$({\bf AFF}_{X})$  Il existe un ensemble fini $T=T(X)$ de places de $k$
tel que, pour tout ensemble fini $S$ de places ne rencontrant pas $T$,
l'image de l'application diagonale $X(k) \to  \prod_{v \in S} X(k_{v})$ est dense.
 
 \medskip

\`A tout corps $k$, \`a toute vari\'et\'e $X$ sur un corps $k$, et plus g\'en\'eralement \`a tout
sch\'ema $X$, on associe son groupe de Brauer--Grothendieck $\Br(X)$ \cite{CTSk21}.

Pour $k$ un corps global, la th\'eorie du corps de classes 
donne des plongements $j_{v} : \Br(k_{v}) \hookrightarrow \Q/\Z$, et une suite exacte fondamentale
$$ 0 \to \Br(k) \to \oplus_{v} \Br(k_{v}) \to \Q/\Z \to 0$$
qui g\'en\'eralise la loi de r\'eciprocit\'e quadratique de Gau{\ss}.

Etant donn\'ee une vari\'et\'e $X$ sur un corps global $k$, 
en utilisant la fonctorialit\'e du groupe
de Brauer, les applications  $j_{v}: \Br(k_{v})  \hookrightarrow  \Q/\Z$ 
induisent un accouplement
$$ X(\A_{k}) \times \Br(X) \to \Q/\Z$$
envoyant un couple $ (\{M_{v}\}, \alpha)$ sur $\sum_{v} j_{v}(\alpha(M_{v}))$.
On note $X(\A_{k})^{\Br} \subset  X(\A_{k})$ le noyau \`a gauche de cet accouplement.
Comme remarqu\'e par Manin en 1970, l'application diagonale $X(k) \to X(\A_{k})$
induit une inclusion
$$ X(k) \subset X(\A_{k})^{\Br}.$$

\medskip

Consid\'erons la propri\'et\'e :

\medskip

$({\bf  BMPH}_{X})$ Soit on a $X(\A_{k})^{\Br}=\varnothing$, soit on a $X(k) \neq \varnothing$.

\medskip

 On dit que l'obstruction de Brauer-Manin au principe de Hasse  
est la seule  pour une classe  
de vari\'et\'es alg\'ebriques d\'efinies sur  $k$  si, pour toute vari\'et\'e $X$ dans cette classe,
on a la propri\'et\'e ${\bf  BMPH}_{X}$.

\medskip

Pour $X/k$ projective, on dit que l'obstruction de Brauer-Manin 
\`a l'approximation faible  
est la seule  pour $X$ si l'on a :

\medskip

$({\bf BMAF}_{X})$  L'ensemble $X(k)$ est dense dans $X(\A_{k})^{\Br}$.

\medskip

Cette propri\'et\'e implique ${\bf BMPH}_{X}$.
Il y a une variante o\`u, pour $k_{v} = \R$ et $k_{v}=\C$, on remplace
$X(k_{v})$ par l'ensemble de ses composantes connexes.

\medskip

Pour les vari\'et\'es projectives, lisses, g\'eom\'etriquement int\`egres
qui sont rationnellement connexes, en particulier celles qui sont g\'eom\'etriquement
unirationnelles, la  propri\'et\'e  ${\bf BMAF}_{X}$ implique la propri\'et\'e d'approximation faible faible  ${\bf AFF}_{X}$.
Ceci r\'esulte du fait que dans ce cas le quotient $\Br(X)/\Br(k)$ est fini.

\medskip

Pour les $k$-vari\'et\'es projectives et lisses g\'eom\'e\-tri\-quement
int\`egres, chacune des propri\'et\'es d\'efinies ci-dessus ne d\'epend que du corps des fonctions de $X$ :
si $X$ et $Y$ sont deux telles $k$-vari\'et\'es birationnellement \'equivalentes, l'une des propri\'et\'es
vaut pour $X$ si et seulement si elle vaut pour $Y$.

\medskip

Si $X$ et $Y$ sont deux $k$-vari\'et\'es, et $Z=X \times_{k}Y$, on a  
$Z(k)= X(k) \times Y(k)$.

Si $f: X \to Y$ est un $k$-morphisme,
il induit une application $X(\A_{k})^{\Br}  \to Y(\A_{k})^{\Br}$.

Si  $X,Y$ sont deux  vari\'et\'es 
projectives et lisses g\'eom\'e\-tri\-quement
int\`egres sur un corps de nombres $k$, et $Z=X \times_{k}Y$, c'est un r\'esultat de Skorobogatov et Zarhin que l'on a
$$Z(\A_{k})^{\Br} = X(\A_{k})^{\Br} \times Y(\A_{k})^{\Br}.$$

\section{Points rationnels des vari\'et\'es rationnellement connexes sur un corps global}\label{pointsrat}

  La conjecture suivante fut faite par Sansuc et moi en 1979  pour les surfaces g\'eom\'etriquement rationnelles \cite{CTSa80},
et \'etendue  aux vari\'et\'es rationnellement connexes en toute dimension en 1999 (voir \cite{CT03}).

\medskip

\begin{conj}\label{conjA}
L'obstruction de Brauer-Manin au principe de Hasse et \`a l'approximation faible
pour les points rationnels est la seule obstruction
pour les vari\'et\'es projectives, lisses, rationnellement connexes
sur un corps global.
\end{conj}

Avec les notations ci-dessus, ceci dit que pour toute vari\'et\'e $X$ projective, 
lisse, rationnellement connexe sur un corps global, on a
 ${\bf  BMPH}_{X}$ et  ${\bf BMAF}_{X}$.

\subsection{\bf Espaces homog\`enes de groupes alg\'ebriques lin\'eaires connexes}

Soit $k$ un corps de nombres. 

\medskip

Pour  $G$ un $k$-groupe semisimple simplement connexe,
  $E$ un espace principal homog\`ene de $G$,  et $X$ une $k$-compactification lisse de $E$,
  des travaux de Eichler, Kneser, Harder et Tchernousov
\'etablirent ${\bf PH}_{X}$. Leurs travaux, et ceux de Platonov, \'etablirent
${\bf  AF}_{X}$. Leurs travaux \'etablissent aussi ${\bf PH}_{X}$ et 
${\bf  AF}_{X}$ pour les vari\'et\'es projectives  $X$ qui sont espaces homog\`enes
d'un groupe alg\'ebrique lin\'eaire connexe~$G$.

Pour $G$ un $k$-groupe alg\'ebrique lin\'eaire connexe quelconque,
$E$ un espace principal homog\`ene de $G$, et $X$ une $k$-compactification lisse  de $E$,
on a les propri\'et\'es ${\bf BMPH}_{X}$ et ${\bf  BMAF}_{X}$ (Voskresenski\v{\i} pour les tores, Sansuc en g\'en\'eral), et donc aussi l'approximation faible faible ${\bf AFF}_{X}$.
Ceci vaut aussi si $E$ est un espace   homog\`ene de $G$ lin\'eaire connexe
lorsque les stabilisateurs g\'eom\'etriques sont connexes (Borovoi).

Par contre, la question suivante est en g\'en\'eral ouverte.

\begin{prob}
 Soit $G$ un $k$-groupe lin\'eaire connexe, $H \subset G$ un $k$-sous-groupe fini.
Soit $X$ une $k$-compactification lisse du quotient $G/H$. A-t-on la propri\'et\'e ${\bf  BMAF}_{X}$,
ou du moins la propri\'et\'e ${\bf AFF}_{X} $ ?
\end{prob}

Comme remarqu\'e par T. Ekedahl et moi en 1988 (voir \cite[Chap. 3]{Se92}), une r\'eponse positive
pour  ${\bf AFF}_{X} $,
appliqu\'ee \`a un groupe fini  abstrait 
$H$ plong\'e dans $GL_{n,k}$ pour $n$ entier
convenable, implique 
que le groupe fini $H$ est le groupe de Galois d'une extension galoisienne
finie de corps $K/k$, propri\'et\'e  qu'on ne sait pas \'etablir pour tous les groupes finis.

Un progr\`es r\'ecent dans cette direction a \'et\'e accompli par Harpaz et Wittenberg  \cite{HW20} pour
une classe de groupes finis $H$ comprenant les groupes nilpotents (constants), ce qui leur permet,
pour ces groupes, de retrouver et pr\'eciser, du point de vue du comportement local,  le th\'eor\`eme de Shafarevich que  les groupes finis
r\'esolubles sont des groupes de Galois sur tout corps de nombres   (on trouve la d\'emonstration
de ce th\'eor\`eme de Shafarevich dans des ouvrages de Ishkhanov, Lur'e, Faddeev et de
Neukirch, Schmidt, Wingberg).

\subsection{\bf Surfaces de del Pezzo et vari\'et\'es de Fano}
Les surfaces de del Pezzo sont les vari\'et\'es de Fano  de dimension 2.

\begin{prob}\label{BMPintcomp} Soit $k$ un corps global.
 Soit $X \subset \P^n_{k}$ une intersection compl\`ete   lisse d\'efinie par l'annulation simultan\'ee
 de  formes homog\`enes  $(f_{1}, \dots, f_{r})$ de degr\'es respectifs $(d_{1}, \dots, d_{r})$.
Si $X$ est dimension au moins 3 et l'on a $d_{1} + \dots + d_{r} \leqslant n$,
a-t-on ${\bf PH}_{X}$ ?
A-t-on ${\bf  AF}_{X}$ ?
\end{prob}

Soit $k$ un corps de nombres.
La m\'ethode du cercle permet d'\'etablir de tels \'enonc\'es pour
$n $ grand  par rapport \`a la somme des $d_{i}$
 (Birch 1961, Schmidt 1985, Skinner 1997).

C'est une exp\'erience commune que pour les vari\'et\'es de Fano il est difficile
d'exhiber des contre-exemples au principe de Hasse.  On consultera
\cite{B18} pour un rapport sur cette direction de recherche tr\`es active.
Sur $k=\Q$, Browning, Le Boudec et Sawin  \cite{BLBS20} ont r\'ecemment montr\'e
que,  si  l'on ordonne (toutes) les hypersurfaces  lisses $X \subset \P^n_{\Q}$ avec $d \leqslant n$ et  $n \geqslant 4$
 par la hauteur des coefficients, alors
100 \% d'entre elles satisfont le principe de Hasse,
 et une proportion positive a des points rationnels.
Pour d'autres r\'esultats ``statistiques'', on consultera \cite{BBL16, LS16, L18, Bri18, SkSo20}.

Consid\'erons maintenant des cas particuliers du probl\`eme \ref{BMPintcomp}.

\begin{prob} Soit $n  \geqslant 4$ et  $X \subset \P^n_{k}$ une hypersurface cubique lisse
sur un corps de nombres. A-t-on ${\bf PH}_{X}$ ?
A-t-on ${\bf  AF}_{X}$ ?
\end{prob}

Si $X$ contient une droite rationnelle $\P^1_{k} \subset \P^n_{k}$,
alors  on a ${\bf  AF}_{X}$   \cite{Har94}.

Pour $k=\Q$,  et $n  \geqslant 9$,  Heath-Brown utilisa la m\'ethode du cercle pour
\'etablir $X(\Q) \neq \varnothing$ pour toute hypersurface cubique lisse,
et C. Hooley \'etablit le principe de Hasse  pour $X$ lisse dans le cas $n=8$.

Sur un corps global  $k$ de caract\'eristique $p>5$, pour $n  \geqslant  5$,
Zhiyu Tian \cite{T17}  a \'etabli ${\bf PH}_{X}$.

 \medskip

Pour $X \subset \P^n_{k}$  intersection compl\`ete lisse sur un corps de caract\'eristique z\'ero,  on a  $\Br(X)/\Br(k)=0$  si $X$ est  de dimension au moins 3.
Dans ce cas,  sur un corps de nombres, on a donc $X(\A_{k})^{\Br}= X(\A_{k})$.
En dimension 2, par exemple pour les surfaces cubiques et les intersections de deux quadriques dans $\P^4_{k}$,
ce n'est plus n\'ecessairement le cas, il faut tenir compte de l'obstruction de Brauer--Manin.

\medskip

\begin{prob} Soit $X \subset \P^3_{k}$ une  surface cubique lisse
sur un corps de nombres. A-t-on ${\bf BMPH}_{X}$ ?
A-t-on ${\bf  BMAF}_{X}$ ?
\end{prob}

 Le cas des surfaces  diagonales  $X \subset \P^3_{\Q}$, d'\'equation
 $ax^3+by^3+cz^3+dt^3=0$, avec $a,b,c,d$ entiers non nuls, 
 sans facteur cubique,  et premiers entre eux dans leur ensemble, a \'et\'e test\'e.  On sait (Cassels-Guy 1966) 
 que ${\bf PH}_{X}$ ne vaut pas toujours pour ces surfaces, 
 mais  dans \cite{CTKS87} on montra que ${\bf BMPH}_{X}$
vaut  lorsque les coefficients sont de valeur absolue plus petite que $100$.
On  a  un  r\'esultat conditionnel, d\^{u} \`a Swinnerton-Dyer \cite{SD01}.
  On suppose la finitude des groupes de Tate-Shafarevich
des courbes elliptiques sur les corps de nombres.
S'il existe un nombre premier $p\neq 3$ qui divise $a$
mais pas $bcd$, et un nombre premier $q \neq 3$ qui divise $b$ mais pas $acd$, alors
le principe de Hasse vaut pour $X$, et ce r\'esultat  conditionnel  implique ${\bf PH}_{X}$ 
pour toute hypersurface cubique diagonale  $X \subset \P^n_{\Q}$ pour $n  \geqslant 4$.

\begin{prob} Soit   $X \subset \P^4_{k}$ une intersection compl\`ete lisse de deux
quadriques sur un corps de nombres $k$.
A-t-on ${\bf BMPH}_{X}$ ?
\end{prob}

On sait que cela vaut si $X$ contient une conique  \cite{Sal88, CT90}.

Par ailleurs, sous l'hypoth\`ese $X(k) \neq \varnothing$, on a
${\bf BMAF}_{X}$  \cite{SaSk91}.

\begin{prob} Soit  $n \geqslant 5$. Soit  $X \subset \P^n_{k}$ une intersection compl\`ete lisse de deux
quadriques sur un corps de nombres $k$. A-t-on ${\bf PH}_{X}$ ?
\end{prob}

 Il est facile de montrer que sous l'hypoth\`ese $X(k)\neq \varnothing$, on a ${\bf{AF}}_{X}$  \cite{CTSaSD87}.
 
On sait que ${\bf PH}_{X}$ vaut  si $X$ contient un ensemble de deux droites conjugu\'ees \cite{CTSaSD87}
ou si $X$ contient une conique (Salberger, 1993, non publi\'e).

On sait que ${\bf PH}_{X}$  vaut pour $n  \geqslant 8$ \cite{CTSaSD87}
  et $n=7$ \cite{HB18}.

Sur un corps global de caract\'eristique $p>2$, ${\bf PH}_{X}$ \'et\'e \'etabli  pour $n \geqslant  5$ par des
m\'ethodes g\'eom\'etriques de d\'eformation par Zhiyu Tian \cite{T17}.

Sur tout corps de nombres, Wittenberg \cite{W07} a donn\'e une preuve conditionnelle
de ${\bf PH}_{X}$ pour $n  \geqslant 5$. Voir    la section \ref{audela} ci-dessous.

\subsection{\bf Espaces totaux de fibrations en vari\'et\'es rationnellement connexes au-dessus de la droite projective}
C'est une classe naturelle de vari\'et\'es \`a consid\'erer si l'on veut \'etablir les r\'esultats par r\'ecurrence
sur la dimension.

\begin{prob}\label{BMenfamille}
Soit $X$ une vari\'et\'e projective et lisse sur un corps de nombres $k$,
munie d'un morphisme $X \to \P^1_{k}$ dont la fibre g\'en\'erique est rationnellement connexe,
et dont les fibres  lisses  $X_{m}$
au-dessus des  $k$-points $m \in \P^1(k)$ satisfont ${\bf BMHP}_{X_{\rm m}}$,
 resp.  ${\bf BMAF}_{X_{\rm m}}$.
A-t-on ${\bf BMHP}_{X} $, resp.  ${\bf BMAF}_{X}$ ?
\end{prob}

Depuis \cite{CTSaSD87}, ce th\`eme a \'et\'e beaucoup explor\'e : travaux de
Skorobogatov, Harari,   Wittenberg,  Harpaz, et de nombreux autres auteurs.
Je renvoie \`a \cite{W18} pour des r\'ef\'erences d\'etaill\'ees.

L'hypoth\`ese sur les fibres est par exemple satisfaite si la fibre
g\'en\'erique $X_{\eta}$ sur le corps $K=k(\P^1)$ est une compactification lisse
d'un espace  homog\`ene d'un $K$-groupe lin\'eaire connexe,
\`a stabilisateurs g\'eom\'etriques connexes (Borovoi).

 \`A une telle fibration on associe 
une mesure de sa complexit\'e arithm\'etique : la somme $\rho$  des degr\'es $[k(m):k] $   des points
ferm\'es $m \in \P^1_{k}$ dont la fibre $X_{m}/k(m)$ est non lisse et ne contient pas de composante
g\'eom\'etriquement int\`egre  de multi\-plicit\'e~1.

On a  des r\'eponses positives  inconditionnelles au probl\`eme \ref{BMenfamille} lorsque $\rho$ est  (tr\`es) petit.
Le  meilleur r\'esultat g\'en\'eral r\'ecent est $\rho \leqslant 3$ \cite{HWW21}.
Pour $\rho$ quelconque, on  a  une r\'eponse conditionnelle positive  \cite{HW16, HWW21}  
si l'on accepte  une conjecture difficile du type de l'hypoth\`ese de Schinzel.
Cette hypoth\`ese,  aussi consid\'er\'ee par   Bouniakovsky,
Dickson,  Hardy et Littlewood, Bateman et Horn, affirme que, pour toute famille finie $P_{i}(t) \in \Z[t]$
de polyn\^{o}mes irr\'eductibles,  \`a coefficients dominants positifs, tels qu'aucun nombre premier ne divise $\prod_{i} P_{i}(m)$ pour tout entier $m$,
  il existe une infinit\'e d'entiers $n$
tels que chaque $P_{i}(n)$ soit un nombre premier.
L'id\'ee d'utiliser l'hypoth\`ese de Schinzel dans ce cadre remonte \`a 1979, et a \'et\'e poursuivie dans divers articles.
Elle vient de conna\^{\i}tre un rebondissement statistique  ``inconditionnel''
 \cite{SkSo20}.

 Un cas simple est donn\'e par une famille de coniques,
d'\'equation affine $$ y^2 - a(t) z^2 -  b(t)=0,$$
 avec $a(t) $ et $b(t)$  polyn\^{o}mes de
degr\'es quelconques. Les fibres de la projection sur l'axe
des~$t$ satisfont le principe de Hasse. Ici $\rho \leqslant 5$
convient.

Depuis  \cite{CTHaSk03} on 
a  aussi beaucoup  \'etudi\'e les \'equations du type
$${\rm Norm}_{K/k}(\Xi) =P (t)$$
avec $\Xi$  ``variable'' dans une extension finie $K/k$ et $P(t) \in k[t]$
polyn\^{o}me non nul. Pour $K/k$ quelconque, les fibres ne satisfont pas en g\'en\'eral
le principe de Hasse mais  elles satisfont  la variante avec obstruction
de Brauer-Manin.

Sur $k=\Q$, des progr\`es fondamentaux en combinatoire additive 
(Green, Tao, Ziegler;  Mathiesen)  ont permis
d'obtenir des r\'esultats inconditionnels avec $\rho$ quelconque.
Les r\'esultats de Green, Tao, Ziegler 
 donnent une version de l'hypoth\`ese de Schinzel pour une famille finie
de formes lin\'eaires \`a deux variables sur $\Q$.
Pour les exemples de vari\'et\'es ci-dessus, \cite{BMSk14, HW16}
montrent ainsi que lorsque $k=\Q$ et que, dans les \'equations ci-dessus, le polyn\^ome
$a(t)b(t)$, resp. le polyn\^{o}me $P(t)$,  a toutes ses racines dans $\Q$, alors on a ${\bf BMAF}_{X} $
(o\`u $X$ d\'esigne un mod\`ele projectif et lisse des vari\'et\'es consid\'er\'ees).

\subsection{Au-del\`a des vari\'et\'es rationnellement connexes}\label{audela}

\medskip

Soit $k$ un corps de nom\-bres.
On ne saurait \'etendre la conjecture \ref{conjA}   \`a toutes les vari\'et\'es projectives et lisses sur $k$,
comme ce fut montr\'e inconditionnellement par Skorobogatov en 1999.
D'autres contre-exemples g\'eom\'etriquement plus simples ont depuis \'et\'e donn\'es.
Cependant, pour  $X$  espace principal homog\`ene   d'une  vari\'et\'e ab\'elienne $A$,
si l'on ignore la composante connexe de l'\'el\'ement neutre aux places archim\'ediennes,
${\bf BMAF}_{X}$
  r\'esulte de la finitude conjecturelle du groupe  de Tate-Shafarevich de la vari\'et\'e ab\'elienne $A$.

  Skorobogatov (2001) conjecture ${\bf BMAF}_{X}$ pour  toute  surface $X$ de type  $K3$.
Si l'on est pr\^et  \`a utiliser non seulement la finitude des groupes de Tate-Shafarevich
mais aussi l'hypoth\`ese de Schinzel,   alors
une m\'ethode sophistiqu\'ee initi\'ee par Swinnerton-Dyer en 1993
permet de pr\'edire un \'enonc\'e de type ${\bf BMHP}_{X}$ pour certaines surfaces $X$ fibr\'ees en courbes
de genre 1 au-dessus de la droite projective. Parmi ces surfaces, on trouve des 
surfaces birationnelles \`a des intersections lisses de deux quadriques dans $\P^4$,
mais aussi des surfaces $K3$. La m\'ethode fut d\'evelopp\'ee dans
\cite{CTSkSD98, W07}. Sous les dites conjectures, Wittenberg  \cite{W07} \'etablit ainsi  ${\bf PH}_{X}$
pour toute intersection compl\`ete lisse $X \subset \P^n_{k}$ pour $n  \geqslant 5$.

La m\'ethode de \cite{SD01}, qui n'utilise ``que'' l'hypoth\`ese de finitude
des groupes de Tate-Shafarevich, a \'et\'e appliqu\'ee par Skorobogatov et Swinnerton-Dyer,
et aussi par Harpaz et Skorobogatov \cite{HS16}, pour \'etudier le principe de Hasse 
pour certaines surfaces de Kummer.

\section{Z\'ero-cycles des vari\'et\'es sur un corps global}

Soit $X$ une vari\'et\'e alg\'ebrique sur un corps $k$. 
L'indice $I(X)$ de la $k$-vari\'et\'e $X$ est  par d\'efinition   le pgcd
des degr\'es $[k(P):k]$ pour tous les points ferm\'es $P$.
C'est aussi le pgcd des degr\'es des extensions finies $K/k$
telles que $X(K) \neq \varnothing$.
Une question plus faible que l'existence d'un point rationnel sur $X$
est celle si l'indice $I(X)=1$.

Dans le cas  des courbes projectives, lisses, g\'eom\'etriquement int\`egres
de genre 0 ou 1,   des quadriques de dimension quelconque,  
et des intersections de deux quadriques,
ces deux questions co\"{\i}ncident, mais ce n'est pas le cas en g\'en\'eral.

Le groupe $Z_{0}(X)$ des z\'ero-cycles sur $X$
est le groupe ab\'elien libre sur les points ferm\'es de $X$.
\`A un z\'ero-cycle $z=\sum_{P} n_{P}P$ ($n_{P} \in \Z$) sur la $k$-vari\'et\'e $X$  on associe
son degr\'e ${\rm deg}_{k}(z):= \sum_{P}n_{P}[k(P):k] \in \Z$.
L'indice $I(X)$ est donc le g\'en\'erateur positif de l'image de l'application
${\rm deg}_{k} : Z_{0}(X) \to \Z$.

Sur un corps de nombres $k$, il est alors naturel de poser la question du principe de Hasse
pour la propri\'et\'e $I(X)=1$ : \'etant donn\'ee une $k$-vari\'et\'e projective, lisse, g\'eom\'etriquement
int\`egre $X$, si on a $I(X_{k_{v}})=1$ pour chaque place $v$, a-t-on alors $I(X)=1$ ?
La r\'eponse est non en g\'en\'eral (courbes de genre 1, intersections compl\`etes  lisses de deux quadriques
dans $\P^4$).

Pour une $k$-vari\'et\'e $X$, on consid\`ere l'accouplement bilin\'eaire
$$Z_{0}(X) \times \Br(X) \to \Br(k)$$
$$ (\sum_{P} n_{P} P,  \alpha ) \mapsto \sum_{P} n_{P} {\rm Cores}_{k(P)/k} (\alpha(P)).$$
Ici $\alpha(P) \in \Br(k(P))$ est l'\'evaluation de $\alpha$ en $P$, et on applique
ensuite la norme, ou corestriction : $\Br(k(P)) \to \Br(k)$.

On peut dans ce cadre d\'efinir une obstruction de Brauer-Manin \`a l'existence d'un z\'ero-cycle de degr\'e 1.
Comme on verra ci-dessous, on peut aussi d\'efinir un analogue de l'obstruction \`a l'approximation faible.

\bigskip

Pour les z\'ero-cycles, on a deux conjectures qui, \`a la diff\'erence de la conjecture  \ref{conjA}, portent sur
{\it toutes} les vari\'et\'es projectives et lisses, sans restriction sur leur g\'eom\'etrie.
Ces conjectures furent faites par Sansuc  et moi (1981) dans le cadre des surfaces rationnelles,
et  \'etendues au cas g\'en\'eral  sous la forme ci-dessous dans  \cite{CT95, CT99}. 
Une conjecture proche mais d'aspect  assez diff\'erent avait
\'et\'e formul\'ee par K. Kato et S. Saito (1983). Voir   \cite{W12}.

\begin{conj}\label{conjB}  Soient $k$ un corps global  et $X$ une $k$-vari\'et\'e projective, lisse,
g\'eom\'etriquement int\`egre sur $k$.
S'il existe une famille de z\'ero-cycles de degr\'e 1 $z_{v} \in Z_{0}(X_{k_{v}})$ tels
que pour tout $\alpha \in \Br(X)$ on ait
$$ \sum_{v} j_{v}(z_{v},\alpha) =0 \in \Q/Z,$$
alors il existe un z\'ero-cycle de degr\'e 1 sur $X$.
\end{conj}

La conjecture est ouverte d\'ej\`a dans le cas des surfaces cubiques lisses.

Soient $X$  et $Y$ des  $k$-vari\'et\'es projectives. Soit $\pi : Y \to X$ un
$k$-morphisme. On lui associe un homomorphisme
$\pi_{*} : Z_{0}(Y) \to Z_{0}(X)$. On dit qu'un z\'ero-cycle sur 
$X$ est rationnellement \'equivalent \`a z\'ero s'il est dans le
sous-groupe de $Z_{0}(X)$ engendr\'e par les $\pi_{*} ({\rm div}_{Y}(g))$
pour $ Y$  variant parmi les courbes normales projectives, la fonction
$g \in k(Y)^{\times}$ variant parmi les fonctions rationnelles non nulles sur une
telle courbe $Y$, et $\pi : Y \to X$ les $k$-morphismes.
Le groupe de Chow $CH_{0}(X)$ des z\'ero-cycles de degr\'e z\'ero sur $X$
est le quotient de $Z_{0}(X)$ par le sous-groupe des z\'ero-cycles rationnellement
\'equivalents \`a z\'ero. Il est muni d'une fl\`eche degr\'e $CH_{0}(X) \to \Z$, dont le noyau
est not\'e $A_{0}(X)$.
L'accouplement $Z_{0}(X) \times \Br(X) \to \Br(k)$ passe au quotient 
par l'\'equivalence rationnelle et induit un accouplement bilin\'eaire
$CH_{0}(X) \times \Br(X) \to \Br(k).$

\medskip

La conjecture suivante englobe la conjecture \ref{conjB}.

\begin{conj}\label{conjC}
 Soient $k$ un corps global  et $X$ une $k$-vari\'et\'e projective, lisse,
g\'eom\'etriquement int\`egre sur $k$.  Le complexe
$$\projlim_{n} CH_{0}(X)/n  \to \prod_{v} \projlim_{n}  CH_{0}(X_{k_{v}})^*/n \to {\rm Hom}(\Br(X), \Q/\Z)$$
induit par la somme des accouplements
du groupe de Brauer de $X$ avec les groupes $CH_{0}(X_{k_{v}}  )$, \`a valeurs
dans $\Br(k_{v}) \subset \Q/\Z$  est une suite exacte.
\end{conj}

On note  $CH_{0}(X_{k_{v}})^*=CH_{0}(X_{k_{v}})$ si $v$ est une place non archim\'edienne,
 puis  $CH_{0}(X_{k_{v}})^*=0$ si $v$ est une place complexe, et pour $v$ r\'eel le quotient de $CH_{0}(X_{\R})$
 par l'image de la norme $CH_{0}(X_{\C}) \to CH_{0}(X_{\R})$.
 
 \medskip

 Le th\'eor\`eme suivant \cite{HW16} est l'aboutissement de  travaux de Salberger \cite{Sal88}, Colliot-Th\'el\`ene, Swinnerton-Dyer, Skorobogatov, Harari \cite{Har94}, Wittenberg \cite{W12}. Un argument relativement \'el\'ementaire mais essentiel du travail  de Salberger \cite{Sal88} avait \'et\'e r\'einterpr\'et\'e par Colliot-Th\'el\`ene et Swinnerton-Dyer (1994) comme une variante inconditionnelle, adapt\'ee aux z\'ero-cycles,  de l'hypoth\`ese de Schinzel mentionn\'ee au paragraphe \ref{pointsrat}.

 \begin{thm} 
  {\rm (Harpaz et Wittenberg)}   Soient $k$ un corps de nombres, $X$
 une $k$-vari\'et\'e projective et lisse g\'eom\'etriquement int\`egre, et
 $\pi : X \to \P^1_{k}$ un $k$-morphisme plat \`a fibre g\'en\'erique
 une vari\'et\'e rationnellement connexe. Si les fibres lisses $X_{m}$
 au-dessus d'un point ferm\'e $m$ de $\P^1$ satisfont la conjecture  \ref{conjB},
 resp. la conjecture \ref{conjC}, alors il en est de m\^eme de $X$.
 \end{thm}

 Y. Liang \cite{Lia13}  avait montr\'e comment on peut \'etablir les conjectures
   \ref{conjB}  et  \ref{conjC}
pour les vari\'et\'es projectives et lisses birationnelles \`a un espace homog\`ene
d'un groupe alg\'ebrique lin\'eaire connexe \`a stabilisateurs connexes
\`a partir du r\'esultat pour les points rationnels (connu gr\^ace \`a Sansuc et Borovoi).

Le th\'eor\`eme ci-dessus  s'applique donc \`a tout $X \to \P^1_{k}$ comme ci-dessus
dont la fibre g\'en\'erique est birationnelle \`a un   espace homog\`ene
d'un groupe alg\'ebrique lin\'eaire connexe \`a stabilisateurs connexes.

 \bigskip
 
 Si $X/k$ est une courbe projective et lisse de genre quelconque, sous  l'hypoth\`ese
 que le groupe de Tate-Shafarevich de la jacobienne $J_{X}$ est fini,
 on a les conjectures    \ref{conjB}  et  \ref{conjC}  pour $X$.

 \medskip
 
Projet ambitieux.  {\it En admettant la finitude des groupes de Tate-Shafarevich des 
vari\'et\'es ab\'eliennes, \'etablir la conjecture    \ref{conjB}   
 pour les surfaces diagonales
$$a x^p + by^p + cz^p+dt^p  = 0$$
de degr\'e $p$ premier dans $\P^3_{\Q}$, par une extension de la m\'ethode
utilis\'ee pour ${\bf PH}_{X}$  et $p=3$ par  Swinnerton-Dyer  \cite{SD01}.}

\medskip

De fa\c con plus g\'en\'erale,  pour une vari\'et\'e projective et lisse,
on souhaiterait ramener la conjecture  \ref{conjB}  au cas
des courbes.
Mais cela semble vraiment hors d'atteinte. Une question plus modeste est :
{\it Suffit-il de conna\^{\i}tre la conjecture  \ref{conjB}  pour toutes les vari\'et\'es projectives et lisses 
de dimension 3 
pour l'avoir en dimension sup\'erieure ?}

\section{Rationalit\'e des vari\'et\'es et invariants birationnels}

Soit $X$ une vari\'et\'e sur un corps $k$. On dit que deux points 
$P,Q \in X(k)$ sont $R$-li\'es s'il existe un ouvert $U \subset \P^1_{k}$
et un $k$-morphisme $U \to X$ avec $P, Q \in  f(U(k))$.
La $R$-\'equivalence sur $X(k)$ est la relation d'\'equivalence
engendr\'ee par cette relation.

\'Etant donn\'es un corps $k$  de caract\'eristique z\'ero 
et une $k$-vari\'et\'e projective, lisse, g\'eom\'etriquement connexe,
l'ensemble $X(k)/R$ et le sous-groupe $A_{0}(X) \subset CH_{0}(X)$ form\'e
des classes de z\'ero-cycles de degr\'e z\'ero, sont des
invariants $k$-birationnels des $k$-vari\'et\'es projectives et lisses,
et ils sont r\'eduits \`a un \'el\'ement si  la $k$-vari\'et\'e $X$ est  
stablement $k$-rationnelle, ou plus g\'en\'eralement  facteur direct
birationnel d'un espace projectif. 

Pour toute $k$-vari\'et\'e $X$ projective,  lisse, g\'eom\'etriquement int\`egre  sur un corps $k$
disons de caract\'eristique z\'ero,   $ i  \geqslant 1$ et $j\in \Z$,
 on dispose des groupes de cohomologie non ramifi\'ee  $H^{i}_{nr}(k(X)/k,\Q/\Z(j))$,
 \`a coefficients dans les racines de l'unit\'e tordues $j$ fois.
Ces groupes sont des invariants $k$-birationnels, r\'eduits \`a $H^{i}(k,\Q/\Z(j))$
si $X$ est  stablement  $k$-birationnelle \`a un espace projectif.
 On consultera \cite{CT19} pour un rapport r\'ecent sur ces invariants.
 On a  $H^{2}_{nr}(k(X)/k,\Q/\Z(1)) = \Br(X)$.
 
 \subsection{Unirationalit\'e}
 
 \begin{prob}
 Soit $k$ un corps infini. Soit $X$ une $k$-vari\'et\'e projective et lisse,
 rationnellement connexe. Supposons $X(k)$ non vide.

(a) L'ensemble $X(k)$ des points rationnels est-il Zariski dense dans $X$ ?

(b) La $k$-vari\'et\'e $X$ est-elle $k$-unirationnelle ?
 \end{prob}
 
 C'est connu   pour les surfaces cubiques lisses, mais ces questions sont ouvertes  pour les surfaces   rationnelles quelconques.
   Pour ces surfaces, une
 r\'eponse affirmative d\'ecoulerait d'une r\'eponse affirmative
 \`a la question suivante \cite[\S V]{CTSa80}:
  \begin{prob}
 Les torseurs universels \cite[\S II.C]{CTSa80} sur les $k$-surfaces rationnelles  projectives et lisses
 sont-ils des $k$-vari\'et\'es (stablement) $k$-rationnelles
 d\`es qu'ils poss\`edent un point rationnel ?
  \end{prob}
 C'est connu pour les surfaces  de Ch\^{a}telet \cite{CTSaSD87}, et plus g\'en\'eralement les sufaces
 fibr\'ees en coniques sur $\P^1_{k}$ avec au plus 4 fibres g\'eom\'etriques singuli\`eres, mais d\'ej\`a le cas des surfaces de del Pezzo $X$ de degr\'e 4  
avec groupe de Picard ${\rm Pic}(X)$ de rang un est ouvert.

 Une r\'eponse affirmative \`a  ce probl\`eme dans le cas $k=\C(\P^1)$  impliquerait
l'unirationalit\'e sur $\C$ des vari\'et\'es complexes de dimension 3 fibr\'ees en coniques
au-dessus du plan projectif $\P^2_{\C}$, ce qui est une question ouverte bien connue.

 \subsection{$R$-\'equivalence}

 Soit $k$ un corps $p$-adique. Soit $f: X \to Y$ un $k$-morphisme
 projectif et lisse de $k$-vari\'et\'es lisses g\'eom\'etriquement int\`egres,
 \`a fibres des vari\'et\'es rationnellement connexes.  
 C'est un th\'eor\`eme de Koll\'ar (1999)
 que pour  $m\in Y(k)$
  l'ensemble $X_{m}(k)/R$ associ\'e \`a la fibre $X_{m}$
  est  fini,
 et que son cardinal est semi-continu sup\'erieurement
 pour la topologie $p$-adique sur $Y(k)$ : pour  $m\in Y(k)$,
 en tout point $n$ d'un  voisinage ouvert convenable de $m$, l'ordre de $X_{n}(k)/R$ 
  est au plus celui de $X_{m}(k)/R$  \cite{K04}.

\begin{prob} {\rm (Koll\'ar)} Sous les conditions ci-dessus, l'ordre de $X_{m}(k)/R$
  est-il localement constant pour la topologie $p$-adique sur $Y(k)$ ?
  \end{prob}

\begin{prob} Soient $k$ un corps parfait de dimension cohomologique 1
  et $X$ une $k$-vari\'et\'e projective, lisse,  (s\'eparablement) rationnellement connexe.
   Supposons $X(k)\neq \varnothing$. 
   
   (a)  L'ensemble  $X(k)/R$  est-il r\'eduit \`a un point  ?
   
   (b) A-t-on  $A_{0}(X)=0$ ?
   
   (c) Ces propri\'et\'es valent-elles au moins si $k$ est un corps $C_{1}$ ?
   \end{prob}

Je renvoie  \`a  \cite[\S 10]{CT11} pour une discussion de divers cas concrets,
tant de corps que de types de vari\'et\'es. La question (b) a une r\'eponse affirmative
 pour les surfaces rationnelles.  La question (c)
a une r\'eponse affirmative pour les intersections compl\`etes lisses de deux quadriques dans $\P^n_{k}$ pour $n  \geqslant 4$. Ici encore, une r\'eponse affirmative \`a  la question (a) dans le cas $k=\C(\P^1)$, et d\'ej\`a  la finitude de
$X(k)/R$ dans ce cas,  impliquerait
l'unirationalit\'e sur $\C$ des vari\'et\'es complexes de dimension 3 fibr\'ees en coniques
au-dessus du plan projectif $\P^2_{\C}$.

  \subsection{Rationalit\'e des intersections de deux quadriques}
  
  Soit $X \subset \P^n_{k}$ une intersection compl\`ete lisse de deux quadriques $f=g=0$ sur
  un corps $k$. Une telle vari\'et\'e est $k$-rationnelle si elle poss\`ede
  une droite $\P^1_{k}$. 
 Un  th\'eor\`eme d'Amer  assure que $X$ contient une droite   
 $\P^1_{k}$ si et seulement si la quadrique  d'\'equation $f+tg=0$ sur le corps $k(t)$
 (o\`u $t$ est une variable)
 contient un $\P^1_{k(t)}$, i.e. si et seulement si
  la forme quadratique $f+tg$ sur le corps  $k(t)$ 
   contient deux hyperboliques.
 
  Pour $k$ alg\'ebriquement clos, on retrouve le fait que $X \subset \P^n_{k}$
  contient une droite si $n  \geqslant 4$, et est rationnelle.
   
   Pour $k$ un corps $C_{1}$, le corps $k(t)$ est $C_{2}$. Dans ce
  cas $X \subset \P^n_{k}$ contient une droite si $n  \geqslant 6$, et est donc $k$-rationnelle.
  C'est le meilleur r\'esultat possible :   pour  $k=\C(z)$ corps des fonctions rationnelles en une variable,
  Hassett et Tschinkel \cite{HT21} donnent un exemple de $X \subset \P^5_{k}$
 qui n'est pas stablement $k$-rationnelle.

   Pour $n=5$, sur un corps quelconque,
un th\'eor\`eme r\'ecent \cite{BW19}, valable sur tout corps,
 dit que la $k$-vari\'et\'e $X$ est $k$-rationnelle si et seulement si elle   contient une droite $\P^1_{k}$.
 On n'a pas par contre de crit\`ere pour la $k$-rationalit\'e stable.
  
  \medskip
  
  Soit $k$ un corps $p$-adique.   Pour  $n  \geqslant 8  $ on   a $X(k) \neq \varnothing$.
Commen\c cons par raffiner certains des r\'esultats de
  \cite[Chap. 3]{CTSaSD87}.
  C'est un th\'eor\`eme  \cite{PS10, HB10, HHK09, L13, PS14}
   que toute   forme quadratique en au moins 9 variables sur un  corps de fonctions d'une variable 
  sur un corps $p$-adique  est isotrope. Ainsi toute forme quadratique en au moins 11 variables
 sur $k(t)$  s'annule sur un vectoriel de dimension 2 sur $k(t)$. 
 Via le  th\'eor\`eme d'Amer, ceci implique
 que pour $n  \geqslant 10$ toute intersection de deux quadriques $X \subset \P^n_{k}$ contient une droite $\P^1_{k}$.
  Donc pour $n  \geqslant 10$, si $X$ est une intersection compl\`ete lisse, elle est $k$-rationnelle.

\begin{prob} Soit $X \subset \P^n_{k}$ une intersection lisse de deux quadriques sur
  un corps $p$-adique $k$.
 Supposons $X(k) \neq \varnothing$. 
  
 (a)  Que peut-on dire sur la $k$-rationalit\'e (stable)  de $X \subset \P^n_{k}$ pour  $6 \leqslant n \leqslant 9$ ?

 (b)  On a $X(k)/R= \{*\}$  pour $n  \geqslant 7$. Que peut-on dire pour $n=5,6$ ?
  
  (c) R\'esultats et questions analogues pour le groupe $A_{0}(X) \subset CH_{0}(X)$ des classes
  de z\'ero-cycles de degr\'e z\'ero.     
    \end{prob}
     
  Si  $p\neq 2$ et  $n=6$, on a $A_{0}(X)=0$ \cite{PS95}.
Par la m\'ethode de sp\'ecialisation  \cite{V15, CTP16} on devrait pouvoir donner des exemples
d'intersections lisses $X$  de deux quadriques dans $\P^5_{k}$, avec $X(k)\neq \varnothing$,
 qui ne sont pas stablement $k$-rationnelles. Voir \cite[Thm. 1.21]{CTP16} et \cite[\S 9, \S 10]{HT21}.
 
  \medskip
  
  On peut aussi se poser la question de la rationalit\'e sur le corps   $\R$ des r\'eels.

\begin{prob}  Soit $X \subset \P^n_{\R}$, $n  \geqslant 4$, une intersection compl\`ete lisse
  de deux quadriques. Supposons que $X(\R)$ est non vide et connexe (pour la topologie r\'eelle).
  Ceci implique-t-il que $X$ est (stablement) $\R$-rationnelle ?
  \end{prob}

  Pour $n=4$, l'hypoth\`ese  implique que $X$ est $\R$-rationnelle. Pour $n=5$,
Hassett et Tschinkel \cite{HT21}  ont montr\'e que $X$ est $\R$-rationnelle si et seulement si
$X$ contient une droite $\P^1_{\R}$. Ainsi pour $n=5$ on peut avoir $X(\R)$ connexe non vide
et $X$ non $\R$-rationnelle. La question de la $\R$-rationalit\'e stable est ouverte.
Pour $n=6$, Hassett, Koll\'ar et Tschinkel (2020) ont montr\'e que $X(\R)$ connexe non
vide \'equivaut \`a  $X$ $\R$-rationnelle. 
En dimensions sup\'erieures le probl\`eme
est ouvert.

  \subsection{Cohomologie non ramifi\'ee}
  
  La question  classique de la rationalit\'e (stable) des hypersurfaces cubiques  lisses dans $\P^4_{\C}$
  am\`ene \`a consid\'erer le probl\`eme suivant,
  qui est li\'e \`a l'\'etude des cycles de codimension deux 
  (voir \cite[Thm. 1.10,  3.1,  3.3]{V15} et \cite[Thm. 5.4, 5.6, 5.8]{CT15}).

  \begin{prob} Soit $X \subset \P^n_{\C}$,  $n  \geqslant 4$ une hypersurface   lisse de degr\'e $d\leqslant n$.
 Soit $K$ un corps contenant $\C$.
Pour $n=4, 5$,  l'application 
 $$H^3(K, \Q/\Z(2)) \to H^3_{nr}(K(X)/K, \Q/\Z(2))$$
 est-elle un isomorphisme pour tout corps $K$ contenant $\C$ ?
 \end{prob}

Pour $n \geqslant 6$ et tout $d\leqslant n$, c'est connu   \cite[Thm. 5.6]{CT15}.
Pour $n=5$ et $d=3$, la r\'eponse est affirmative, cela r\'esulte  \cite[Thm. 5.8]{CT15}  d'un th\'eor\`eme de C.~Voisin (2006).
Pour $n=4$ et $d=3$,  c'est un probl\`eme en g\'en\'eral ouvert \cite{V15}.

\medskip

 Soient $\F$ un corps fini et $\ell$ un nombre premier diff\'erent de la caract\'eristique de $\F$.
  Soit $X/\F$ une vari\'et\'e projective et lisse g\'eom\'etriquement connexe
   de dimension~$d$.
  Le groupe $H^3_{nr}(\F(X)/\F,\Q_{\ell}/\Z_{\ell}(2))$ est un analogue sup\'erieur de la partie $\ell$-primaire
  du groupe de Brauer $\Br(X)$ d'une vari\'et\'e $X/\F$. 
  C'est une extension d'un groupe fini par un groupe divisible.
   Pour $d=2$, i.e. $X$  une surface, on a $H^3_{nr}(\F(X)/\F,\Q_{\ell}/\Z_{\ell}(2))=0$ (corps de classes sup\'erieur).
  A. Pirutka a donn\'e des exemples de vari\'et\'es g\'eom\'etriquement rationnelles de dimension 5    avec
  $H^3_{nr}(\F(X)/\F,\Q_{2}/\Z_{2}(2))\neq 0$. F. Scavia et F. Suzuki viennent de donner un exemple de
  vari\'et\'e de dimension 4 pour laquelle ce groupe est non nul.

 \begin{prob} Pour  toute vari\'et\'e projective et lisse int\`egre $X$ de dimension 3, 
le groupe $H^3_{nr}(\F(X)/\F,\Q_{\ell}/\Z_{\ell}(2))$ est-il divisible ? Est-il nul ?
Est-ce d\'ej\`a le cas pour les vari\'et\'es rationnellement connexes ?
\end{prob}

 On sait l'\'etablir pour quelques classes int\'eressantes de vari\'et\'es :
  les vari\'et\'es fibr\'ees en coniques au-dessus d'une surface \cite{PS16}, 
  et les hypersurfaces cubiques lisses dans $\P^4_{\F}$.
 La question est li\'ee  \`a une forme forte de la conjecture de Tate enti\`ere pour les 1-cycles sur les vari\'et\'es de dimension 3 sur un corps fini, et  \`a la validit\'e de la conjecture  \ref{conjC} 
 ci-dessus pour les surfaces sur un corps global de caract\'eristique positive \cite{CT99, CTK13}.
 Le lien entre la conjecture de Tate enti\`ere pour les 1-cycles sur les vari\'et\'es  sur un corps fini
 et la conjecture  \ref{conjB}  sur un corps global de caract\'eristique positive
 avait \'et\'e fait par S. Saito en 1989.
 On sait \'etablir $H^3_{nr}(\F(X)/\F, \Z/2)= 0$ pour $X \subset \P^4_{\F}$ une hypersurface cubique lisse
 (${\rm car.}(\F) \neq 2$).
 
 \begin{prob} Soit $\F$ un corps fini, ${\rm car}(\F) \neq 2$, 
et soit $X \subset \P^5_{\F}$ une hypersurface cubique lisse.
A-t-on $H^3_{nr}(\F(X)/\F, \Z/2)= 0$ ?
\end{prob}

 	\section{Points rationnels et indice des vari\'et\'es alg\'ebriques.}
	
\begin{prob} Soit $X\subset \P^n_{k}, n  \geqslant 4, $ une   intersection
compl\`ete lisse de deux quadriques sur un 
  corps $k$ de dimension cohomologique 1.
 Pour $n=5$,  a-t-on $X(k) \neq \varnothing$ ?
 \end{prob}

  Pour   $n  \geqslant 6$, c'est vrai et facile.
  Pour $n=4$, la r\'eponse est n\'egative. La
  d\'emonstration repose sur la construction de tr\`es
  grands corps.
 
 \medskip
 
 Soit $C$ une courbe g\'eom\'etriquement int\`egre sur le corps
 des r\'eels  avec $C(\R)=\varnothing$, par exemple la conique 
 d'\'equation homog\`ene  $x^2+y^2+t^2=0$. On sait que le
 corps $\R(C)$ est de dimension cohomologique 1.
 C'est une question ouverte si c'est un corps~$C_{1}$.
 Plus g\'en\'eralement on demande s'il y  a un analogue
 du th\'eor\`eme de Graber, Harris et Starr :
 
 \begin{prob} 
 Toute vari\'et\'e rationnellement connexe $X$  sur le corps $K=\R(C)$
 a-t-elle un point rationnel ?
 \end{prob}

On ne sait d\'ej\`a pas si pour toute vari\'et\'e  projective et lisse rationnellement connexe
$X$ sur $\R$  avec $X(\R) =\varnothing$ il existe un $\R$-morphisme
de la conique sans point vers $X$.

\begin{prob} Existe-t-il un entier $n \geqslant 4$  tel que  
    toute  hypersurface cubique lisse $X \subset \P^n_{k}$  sur un corps $k$ 
  de dimension cohomologique 1  poss\`ede un point rationnel, 
  ou du moins satisfasse $I(X)=1$ ? 
  \end{prob}

Les corps $p$-adiques ne sont pas des corps $C_{2}$.  
On a cependant la question :

\begin{prob}  {\rm (Kato  et Kuzumaki)} 
Pour toute hypersurface $X \subset \P^n_{k}$
de degr\'e $d$  sur un corps $p$-adique, si l'on a  $n\geq d^2$, a-t-on $I(X)=1$ ?
\end{prob}

Ceci a \'et\'e \'etabli par Kato et Kuzumaki (1985)  lorsque le degr\'e $d$ est un
nombre premier. C'est ouvert d\'ej\`a pour $d=4$.

\begin{prob} {\rm (Cassels et Swinnerton-Dyer)} Soient $k$ un corps et $X \subset \P^n_{k}$ une
  hypersurface cubique. Si l'on a $I(X)=1$, a-t-on $X(k) \neq \varnothing$ ?
  \end{prob}

 D. Coray (1976) montra qu'il en est ainsi sur un corps $p$-adique. Pour  
 $X \subset \P^3_{k}$  
 une surface cubique lisse sur un corps quelconque, il montra que l'hypoth\`ese $I(X)=1$ entra\^{\i}ne
  l'existence sur $X$ d'un point ferm\'e de degr\'e 1, 4 ou 10.  La question si on peut \'eliminer $10$ et $4$
   est rest\'ee ouverte. L'analogue de la question pour les surfaces de del Pezzo de degr\'e 2   a une r\'eponse
n\'egative (Koll\'ar et Mella).
 
   \begin{prob} {\rm (Serre)}
   Soient $k$ un corps, $G$ un groupe alg\'ebrique lin\'eaire connexe sur $k$,
   et $E$ un espace principal homog\`ene sous $G$.  Si l'indice $I(E)$ est \'egal \`a $1$,
   a-t-on $E(k)\neq \varnothing$ ?
      \end{prob}
 
 Pour les  espaces homog\`enes non principaux, la propri\'et\'e ne vaut pas,
 des contre-exemples ont \'et\'e construits par Florence et par Parimala.
 
\begin{prob} Soient $k$ un corps, ${\rm car}(k)=0$, et $X \subset \P^4_{k}$ 
 une hypersurface cubique lisse sans point rationnel, i.e.  d'indice $I(X)=3$.
 Soit $Y/k$ une $k$-vari\'et\'e projective lisse g\'eom\'etriquement connexe 
 de dimension au plus 2. S'il existe  une $k$-application rationnelle de $X$ vers $Y$,
 a-t-on $I(Y)=1$ ? \end{prob}

 Par la classification $k$-birationnelle des surfaces g\'eom\'etriquement rationnelles, 
 le probl\`eme se ram\`ene au cas o\`u $Y$ est une surface cubique lisse $k$-minimale.
 Ce cas semble r\'esister aux formules de degr\'e \`a la Rost
  \cite{M03, Z10} qui avaient permis d'\'etendre  le th\'eor\`eme d'Hoffmann  \cite{H95}
  restreignant les dimensions possibles pour les couples de quadriques anisotropes
  admettant une application rationnelle entre elles.

\section{Groupes alg\'ebriques lin\'eaires}

Soit $G$ un groupe alg\'ebrique lin\'eaire r\'eductif connexe sur un corps $k$.
L'ensemble $G(k)/R$ est naturellement muni d'une structure de groupe.
Tout \'el\'ement de ce groupe est d'ordre fini. Si $K/k$ est une extension
transcendante pure, l'homomorphisme $G(k)/R \to G(K)/R$ est  bijectif.

Soit $D$ une alg\`ebre   centrale simple (de rang fini) sur un corps $k$.
Soit $G=SL_{1,D} \subset GL_{1,D}$ le groupe alg\'ebrique des \'el\'ements
de norme r\'eduite 1. C'est un $k$-groupe alg\'ebrique semisimple simplement connexe.
Un th\'eor\`eme de Voskresenski\u{\i} utilisant un th\'eor\`eme de Platonov
identifie dans ce cas $G(k)/R$ au groupe $SK_{1}(D)$ quotient
du groupe  des \'el\'ements de $D^{\times}$ de norme r\'eduite 1 par le groupe
$[D^{\times}, D^{\times}]$ engendr\'e par les commutateurs. 

Des travaux de Platonov, Yanchevski\u{\i}, Merkurjev,    Chernousov
ont identifi\'e  le quotient $G(k)/R$ pour beaucoup de groupes classiques,
tant simplement connexes qu'adjoints, et ont au passage \'etabli 
sa commutativit\'e. Le probl\`eme  g\'en\'eral suivant reste cependant ouvert.

\begin{prob} Soient $k$ un corps et $G$ un $k$-groupe
alg\'ebrique r\'eductif connexe. Le groupe quotient $G(k)/R$
est-il commutatif ?
\end{prob}

Voskresenski\v{\i} (1977) avait pos\'e la question 
pour $G$ lin\'eaire connexe sur un corps quelconque.
Pour un groupe lin\'eaire connexe non r\'eductif, sur un corps
non parfait, F. Scavia (2021) a donn\'e une r\'eponse n\'egative.

\medskip

Pour les deux probl\`emes suivants, on consultera le rapport de P. Gille \cite{Gi07}.

\begin{prob}
Soient $k$ un corps et $G$ un  $k$-groupe r\'eductif connexe.
Si pour tout corps $K$ contenant $k$, le groupe $G(K)/R$
est trivial, ceci implique-t-il que $G$ est facteur direct birationnel
d'une $k$-vari\'et\'e $k$-rationnelle ?
\end{prob}

\begin{prob} Soit $k$ un corps parfait de dimension
cohomologique $\leqslant 3$. Si $G$ est un $k$-groupe
semisimple simplement connexe, a-t-on  $G(k)/R=1$ ?
\end{prob}
Ce probl\`eme  est motiv\'e par les travaux de Suslin
 sur le groupe $SK_{1}(D)$ d'une alg\`ebre simple centrale.
On a un certain nombre de r\'esultats lorsque la dimension
cohomologique est $\leqslant 2$.

 \begin{prob}
Soit $G$ un groupe alg\'ebrique lin\'eaire connexe sur le corps $\R$ des r\'eels.
La $\R$-vari\'et\'e $G$ est-elle $\R$-rationnelle, i.e. le corps des fonctions de $G$
est-il transcendant pur sur $\R$ ?
\end{prob}

\begin{prob}  Soit $k$ un corps de type fini sur $\Q$.
Si $G$ est un $k$-groupe lin\'eaire connexe,  le quotient 
$G(k)/R$ est-il fini ? 
\end{prob}

C'est connu pour $G$ un $k$-tore (Colliot-Th\'el\`ene et Sansuc 1977) et pour
$k$ un corps de nombres (P. Gille 1997).
Pour $G=SL_{1,D}$ le $k$-groupe alg\'ebrique des \'el\'ements
de norme 1 dans une alg\`ebre  centrale  simple $D$ sur $k$,
le probl\`eme se traduit ainsi :

\begin{prob} Si $D$ est une alg\`ebre centrale simple  sur un corps $k$ de type fini sur $\Q$,
le groupe $SK_{1}(D)$ est-il fini ?
\end{prob}

\begin{prob} Soit $H/\C$ un groupe lin\'eaire connexe
et $G \subset H$ un sous-groupe ferm\'e  connexe.
Le quotient $H/G$ est-il une vari\'et\'e rationnelle ?
\end{prob}

C'est une question c\'el\`ebre, d\'ej\`a  pour $H=GL_{n,\C}$
et $G=PGL_{m,\C}$. Pour traiter cette question, on peut essayer 
d'utiliser la cohomologie non ramifi\'ee.

\medskip

Soit $k=\C$, et soit $X$ une compactification lisse de  $GL_{n,\C}/G$, avec $G \subset GL_{n,\C}$
sous-groupe alg\'ebrique ferm\'e connexe.
Pour tout corps $K$ contenant $\C$,  et $i=1,2$  on sait
que  $$H^{i}(K,\Q/\Z(i-1)) =    H^{i}_{nr}(K(X)/K,\Q/\Z(i-1)).$$
Pour $i=2$, ceci dit que le groupe de Brauer   de $X\times_{\C}K$
est r\'eduit \`a l'image de $\Br(K)$, \'enonc\'e essentiellement d\^{u}
\`a Bogomolov.
Dans une s\'erie d'articles, Merkurjev \cite{M17}  et Sanghoon Baek \cite{B21}  ont \'etabli
$H^{3}_{nr}(\C(X)/\C, \Q/\Z(2))=0$ pour de nombreuses classes
de groupes r\'eductifs $G$.

\begin{prob} Dans chacun de ces cas,   pour tout corps $K$
contenant $\C$, la fl\`eche
$H^{3}(K,\Q/\Z(2)) \to   H^{3}_{nr}(K(X)/K,\Q/\Z(2))$
est-elle un isomorphisme ?
\end{prob}
 
Avec les m\'ethodes d\'ecrites dans  \cite[\S 5]{CT15},
on pourrait essayer de r\'esoudre ce probl\`eme via l'\'etude
des cycles de codimension deux  d'une bonne
compactification lisse  de $GL_{n,\C}/G$.

\end{document}